\documentclass{amsart}

\usepackage{geometry} 
\usepackage{graphicx}
\usepackage{adjustbox}
\usepackage{amssymb}
\usepackage{epstopdf}
\usepackage{amsmath,amsthm,epsfig,amssymb}
\usepackage{color}
\usepackage{adjustbox}

\DeclareMathOperator{\lcm}{lcm}

\newtheorem{theorem}{Theorem}[section]

\newtheorem{corollary}[theorem]{Corollary}
\newtheorem{proposition}[theorem]{Proposition}
\newtheorem{conjecture}[theorem]{Conjecture}

\theoremstyle{definition}
\newtheorem{definition}[theorem]{Definition}

\theoremstyle{remark}

\numberwithin{equation}{section}

\begin{document}

\title{A geometric consideration of the Erd\H{o}s-Straus conjecture}

\author{K.~Bradford}
\address{Department of Mathematics and Statistics \\
University of Nevada, Reno \\
1664 N. Virginia St. \\
Reno, NV 89557-0084, USA \\
\ } 
\email{kylebradford@unr.edu}

\author{E.~Ionascu}
\address{Department of Mathematics \\
Columbus State University \\
4225 University Ave. \\
Columbus, GA 31907, USA}
\email{ionascu@columbusstate.edu}

\date{December 7, 2014}

\begin{abstract}
In this paper we will explore the solutions to the diophantine equation in the Erd\H{o}s-Straus conjecture.  For each prime $p$  we are discussing the relationship between the values $x,y,z \in \mathbb{N}$  so that $$ \frac{4}{p} = \frac{1}{x} + \frac{1}{y} + \frac{1}{z}.$$
\ 

\noindent We will separate the types of solutions into two cases.  In particular we will argue that the most common relationship found is  $$ x = \left\lfloor \frac{py}{4y-p} \right\rfloor + 1.$$
\ 

\noindent Finally, we will make a few conjectures to motivate further research in this area.
\end{abstract}

\maketitle

\section{Introduction} \label{introduction}

The Erd\H{o}s-Straus conjecture suggests that for every $n \geq 2$  there there exist natural numbers $x,y,z \in \mathbb{N}$  so that 

\begin{equation} \label{eq: main}
\frac{4}{n} = \frac{1}{x} + \frac{1}{y} + \frac{1}{z}.
\end{equation}
\ 

Naturally this reduces to prime numbers.  This means that a sufficient condition for proving the conjecture is if one could show that for every prime number $p$, there exist natural numbers $x,y,z \in \mathbb{N}$ that satisfy (\ref{eq: main}).  It is safe to assume that $x \leq y \leq z$  as one of the values will be the largest and one will be the smallest.  The solutions to (\ref{eq: main}) need not be unique.  For example we see that 
\begin{equation} \label{eq: seventeen}
\begin{aligned}
\frac{4}{17} &= \frac{1}{5} + \frac{1}{34} + \frac{1}{170} \\
&= \frac{1}{5} + \frac{1}{30} + \frac{1}{510} \\
&= \frac{1}{6} + \frac{1}{15} + \frac{1}{510} \\
&= \frac{1}{6} + \frac{1}{17} + \frac{1}{102}.
\end{aligned}
\end{equation}
\ 

\begin{figure}[h] 
\center

\fbox{\reflectbox{\rotatebox[origin=c]{180}{\adjustbox{trim={.07\width} {.08\height} {0.085\width} {.12\height},clip}{\includegraphics[scale=0.4]{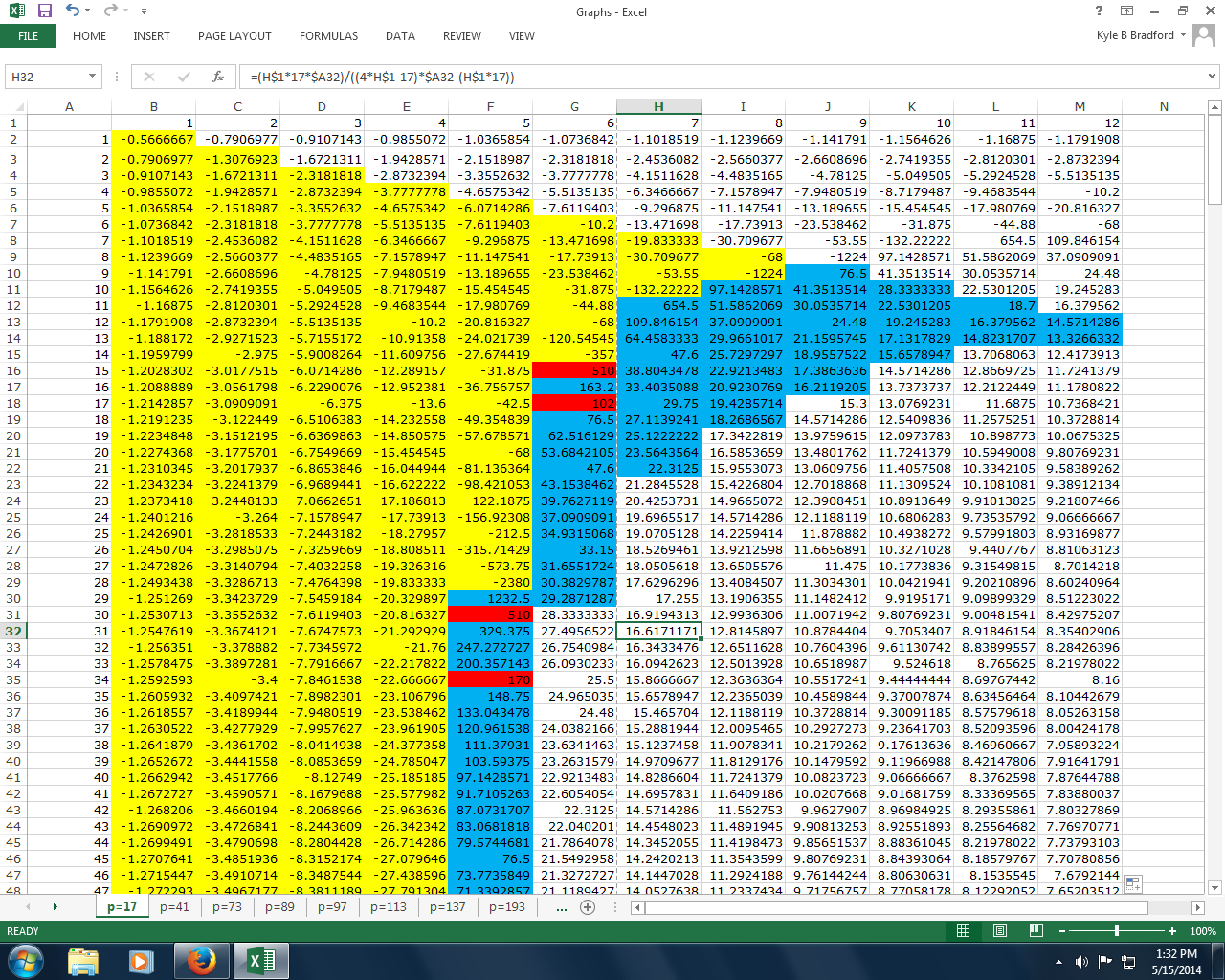}}}}}

\caption{For $p = 17$ the cells represent the standard xy integer lattice with the colored cell entries representing the $z$ values.  If $x > y$  or $y > z$ the cell color is white.  If $z$  is negative the cell color is yellow.  If $z$  is a solution the cell color is pink.  Otherwise the cell color is blue.} \label{fig: seventeen}
\end{figure}

The Erd\H{o}s-Straus conjecture dates back to the 1940s and early 1950s \cite{pe1, ob, ro}.  People have attempted to solve this problem in many different ways.  For example, algebraic geometry techniques to give structure this problem (see \cite{cs}), analytic number theory techniques to find mean and asymptotic results (see \cite{ec, et, dl, san1, san2, v, wweb1, yang}), comparing related fractions, such as $k \slash n$ for $k \geq 2$ (see \cite{aaa, ec, gm, rav, wweb2, wweb3}), computational methods (see \cite{sa}), organizing primes $p$  into two classes based off of the decompositions of $4 \slash p$  in hopes to find a pattern within each class (see \cite{bh, et, san1, san2}), and looking for patterns in the field of fractions of the polynomial ring $\mathbb{Z}[x]$ instead of $\mathbb{Q}$ (see \cite{s}) just to name a few.  The current authors have made attempts to make equivalent conjectures in different number fields \cite{bi}.  The best-known method was developed by Rosati \cite{ro}.  Mordell \cite{mo}  has a great description of this method and many attempts use the techniques in his paper (see \cite{iw, san3, t, y}).  

The purpose of this paper is to classify each solution based on its geometric location.  Figure \ref{fig: seventeen} shows the geometric location of the solutions listed in (\ref{eq: seventeen}) as pink cells where the cells represent the standard xy integer lattice when both $x>0$  and $y>0$.  This image was made with a Microsoft excel worksheet by using conditional formatting of the cell colors.  The pink cells that border the yellow cells in figure \ref{fig: seventeen} will be of particular interest.  In this case we see that all the pink cells border the yellow cells.  To define the border between the yellow and blue cells in figure \ref{fig: seventeen} we need to relate $x$  and $y$.  We let the cells be white if $x > y$  or if $y > z$.  When $p=17$  we will see that $y>z$  if $y > 34x \slash (4x - 17)$.   We will let the cells be yellow if $z < 0$.  The cells will be yellow if $y < 17x \slash (4x - 17)$  and the cells will be blue or pink if $17x \slash (4x-17) \leq y \leq 34x \slash (4x-17)$.  The cells are pink only if $z$  is an integer.   Our main argument will be that a overwhelming majority of the solutions fall along the boundary of all $(x,y)$  values that give $z > 0$.  

We will also see that for all primes $p \neq 2$  and $p \neq 2521$  there exists at least one solution to (\ref{eq: main}) so that $x = \left\lfloor py \slash (4y-p) \right\rfloor +1$, $\gcd \left(p,y \right) = 1$  and $z = p \cdot \lcm \left( x, y \right)$.  For $p = 17$  we see that there are two solutions with this pattern.  These solutions are 

\begin{align*}
\frac{4}{17} & = \frac{1}{5} + \frac{1}{30} + \frac{1}{510} \\
&= \frac{1}{6} + \frac{1}{15} + \frac{1}{510}.
\end{align*}
\ 

Finally we will see that for all primes $p \not\in \{ 2,3,7,47,193,2521\}$  there exists at least one solution to (\ref{eq: main})  so that $y = \left\lfloor px \slash (4x - p) \right\rfloor + 1$, $\gcd(p,y) = 1$  and $z = p \cdot \lcm \left( x,y \right)$.   For $p = 17$  we see that there is only one solution with this patter.  This solution is $$ \frac{4}{17} = \frac{1}{6} + \frac{1}{15} + \frac{1}{510}.$$
\ 

The rest of the paper is organized as follows:  in section \ref{sec: main} we will describe the main results of our paper without proof and in section \ref{sec: development} we will fill in the necessary details.

\section{Main Results} \label{sec: main}
We will generalize the results made in the introduction to any prime $p$.  Our first goal in this endeavor is to define the boundary between the yellow cells and the blue or pink cells as in figure \ref{fig: seventeen} for a general prime $p$.  We notice that if $$y < \frac{px}{4x-p}$$

\noindent then 
\begin{equation} \label{eq: border inequality}
\frac{4}{p} < \frac{1}{x} + \frac{1}{y}.
\end{equation}
\ 

To solve (\ref{eq: main}) when (\ref{eq: border inequality}) holds, we necessarily need $z$  be negative. Because this cannot happen, this implies that $$y \geq \frac{px}{4x-p}.$$

To solve (\ref{eq: main}), the equation $ 4xy - p(x+y) = 0$  cannot hold because if it did hold, then $$ \frac{4}{p} = \frac{1}{x} + \frac{1}{y} $$
\ 

\noindent and necessarily $z$  cannot be an integer.  This equation will, however,  define the boundary between the yellow cells and the blue or pink cells mentioned from figure \ref{fig: seventeen}  and it will apply to any prime $p$.  To be on the correct side of this boundary we see that 

\begin{equation} \label{eq: boundary y inequality}
4xy-p(x+y) > 0.
\end{equation}
\ 

To be along the boundary, yet satisfy (\ref{eq: boundary y inequality}), we need to select the integer values of $x$  and $y$  so that the left hand side of the inequality (\ref{eq: boundary y inequality}) is the smallest possible positive value.  The following definition will describe two ways that a solution to (\ref{eq: main}) can be along this boundary.

\begin{definition} \label{def: type I}

A solution to (\ref{eq: main})  is a type I(a) solution if 
\begin{equation} \label{eq: type I(a)}
y = \left\lfloor \frac{px}{4x-p} \right\rfloor + 1.
\end{equation}
\ 

\noindent A solution to (\ref{eq: main})  is a type I(b) solution if 
\begin{equation} \label{eq: type I(b)}
x = \left\lfloor \frac{py}{4y-p} \right\rfloor + 1.
\end{equation}
\ 
\noindent A solution is called a type I solution if it is a type I(a) solution, a type I(b) solution or both.
\end{definition}

If we relate this to figure \ref{fig: seventeen}, then type I solutions are given by the pink cells that border a yellow cell from the bottom or from the right.  In particular, a type I(a) solution is given by a pink cell that borders a yellow cell from the bottom and a type I(b) solution is given by a pink cell that borders a yellow cell from the right.  We quickly find a relationship between type I(a) solutions and type I(b) solutions, which we outline in the following proposition.

\begin{proposition} \label{prop: three}
If a solution is a type I(a)  solution then it is a type I(b) solution.
\end{proposition}

This means that if a solution to (\ref{eq: main})  is of type I, then it is of type I(b).  We can use the two terms interchangeably.  There is computational evidence to suggest that the only prime $p$  where there is no solution of type I(a) is when $p = 193$.  This computation evidence is through all primes less that $10^{8}$.  We summarize this conclusion in the following conjecture.

\begin{conjecture}
The only prime $p$  where there is no solution of type I(a) is $p = 193$.
\end{conjecture}

Because all type I(a) solutions are type I(b) solutions, we can make a stronger statement about type I(b) solutions.  Because $$ \frac{4}{193} = \frac{1}{50} + \frac{1}{1930} + \frac{1}{4825} $$

\ 

\noindent is a type I(b) solution, there is computational evidence to suggest that every prime $p$  has a solution of type I(b).  This computational evidence is through all primes less that $10^{8}$. We summarize this conclusion in the following conjecture.

\begin{conjecture}
Every prime $p$  has a solution of type I(b).
\end{conjecture}

The fact that every prime has at least one solution of type I(b) gives the authors of this paper the impression that the proof of the Erd\H{o}s-Straus conjecture reduces to to finding a solution of type I(b) for every prime $p$.  This may not be true, but it leads us to ask some natural questions about which primes $p$ have a decomposition that we can prove are of type I(b).  First we recall a theorem from \cite{iw}.

\begin{theorem} \label{thm: ionascu} {\bf Ionascu-Wilson}
Equation (\ref{eq: main})  has at least one solution for every prime number $p$, except possible for those primes of the form $p \equiv r \mod 9240$  where $r$  is of the $34$  entries in the table:

\ 

\center 
\begin{tabular}{| c | c | c | c | c | c |} \hline
$1$ & $169$ & $289$ & $361$ & $529$ & $841$ \\ \hline
$961$ & $1369$ & $1681$ & $1849$ & $2041$ & $2209$ \\ \hline
$2521$ & $2641$ & $2689$ & $2809$ & $3361$ & $3481$ \\ \hline
$3529$ & $3721$ & $4321$ & $4489$ & $5041$ & $5161$ \\ \hline
$5329$ & $5569$ & $6169$ & $6241$ & $6889$ & $7561$ \\ \hline
$7681$ & $7921$ & $8089$ & $8761$ &  &  \\ \hline
\end{tabular}
\end{theorem}

\ 

The decompositions created to prove theorem \ref{thm: ionascu} were given in \cite{iw} and can be tested to determine whether or not they were of type I(b).  The following theorem tells us that every solution provided is of type I(b)

\begin{theorem} \label{thm: two}
Every prime $p$  that is guaranteed a solution by theorem \ref{thm: ionascu} has at least one solution of type I(b).
\end{theorem}

Although every prime has at least one solution of type I(b), we were curious to know whether or not every solution was of type I(b). We can see for $p=17$  that every solution was of type I(b), however, for other primes there exist solutions that are not of type I.  For example, we have that $$ \frac{4}{71} = \frac{1}{20} + \frac{1}{284} + \frac{1}{355} $$
\ 

\noindent where we see that 
\begin{align*}
x &= \left\lfloor \frac{71 \cdot 284}{4 \cdot 284 - 71} \right\rfloor + 2 \\
&= 20.
\end{align*}
\ 

To account for the remaining solutions, we make the following definition.

\begin{definition} \label{def: type II}
A solution to (\ref{eq: main})  that is not a type I solution is a type II solution.
\end{definition}

\begin{figure}[h] 
\center
\fbox{\includegraphics[scale=0.75]{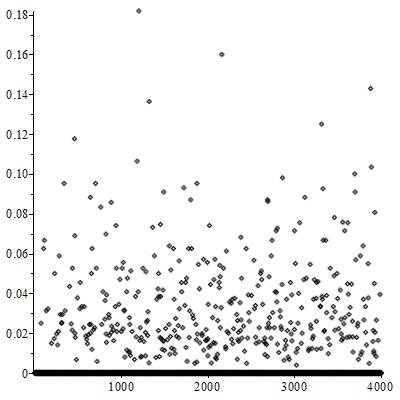}}

\caption{This graph shows the proportion of type II solutions for each prime $p$.} \label{fig: proportion II solutions}
\end{figure}
\ 

It is natural to ask if there is a pattern within the class of type II solutions.  Although there is most likely no upper bound to the number of type II solutions that exist for a given prime, it appears that as the number of type II solutions grow, the number of type I solutions grow as well.  They do not, however, appear to grow at a uniform rate.  Figure \ref{fig: proportion II solutions}  shows the proportion of type II solutions for each prime less than $4000$.  There is no prime less than $4000$  that has less than $80$\%  of its solutions of type I, but this proportion seems sporadic.

\begin{figure}[h]
\center
\fbox{\includegraphics[scale=0.75]{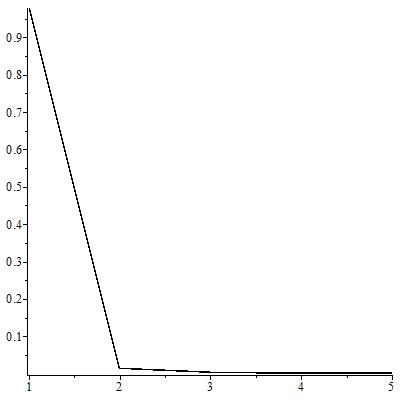}}

\caption{This graph draws the probability distribution function defined from table \ref{tab: empirical}.  The points in the pdf are connected with lines.} \label{fig: empirical}
\end{figure}

We can see from figure \ref{fig: proportion II solutions}  that most primes have no type II solutions at all, so our next goal was to make an empirical distribution for the solutions to (\ref{eq: main})  based on the proximity of the solution to the boundary.  For example, there are $38434$  solutions to (\ref{eq: main})  for primes $p \leq 4000$.  We will separate the number of solutions to (\ref{eq: main})  for prime numbers $p$  into categories based on whether the solutions satisfy $$x = \left\lfloor \frac{py}{4y-p} \right\rfloor + i$$
\ 

\noindent for $1 \leq i \leq 5$.  Table \ref{tab: empirical}  and figure \ref{fig: empirical}  summarize what we have found for primes $p \leq 4000$. \\

This distribution shows our point very well.  If we are to describe a pattern for solutions to (\ref{eq: main}) for a general prime $p$, it appears that it is a safe assumption to let $$ x = \left\lfloor \frac{py}{4y-p} \right\rfloor + 1.$$
\ 

Next we turn our attention to another pattern one can easily identify for solutions of (\ref{eq: main}).  As mentioned in the introduction, we can see that for all primes $p$  such that $p \neq 2$  and $p \neq 2521$  there exists a solution so that $x = \left\lfloor py \slash (4y-p) \right\rfloor + 1$, $\gcd \left(p,y \right) = 1$  and $z = p \cdot \lcm \left( x,y \right)$. This has been checked computationally for all primes less than $10^{8}$.  Instead of trying to explain why the two primes $p=2$  and $p=2521$  do not follow this pattern, we argue that it suffices to find a prime $p^{*}$  large enough so that every prime larger than $p^{*}$  has the pattern we describe above.  This brings up two conjectures.  We believe that these conjectures govern at least one way to find a general pattern for the solutions of (\ref{eq: main}). \\

\begin{table}[h]
\begin{center}
\begin{tabular}{| c | c | c |}
\hline
$i$ & \# solutions & proportion \\ \hline \hline
$1$ & $37612$ & $0.9786$ \\ \hline
$2$ & $517$ & $0.0135$ \\ \hline
$3$ & $170$ & $0.0044$ \\ \hline
$4$ & $64$ & $0.0017$ \\ \hline
$5$ & $71$ & $0.0018$ \\ \hline
\end{tabular} 
\caption{This table shows the empirical probability distribution function of the solutions to (\ref{eq: main}) based on their proximity to the boundary values that make $z$ positive.  The solutions are accumulated for primes less than $4000$ and separated into categories based on whether the solutions satisfy $x = \left\lfloor py \slash (4y-p) \right \rfloor + i$.}  \label{tab: empirical}
\end{center}
\end{table}
\ 

First we mention that for any prime $p \neq 2$  and $y \in \mathbb{N}$  that  satisfy (\ref{eq: main})  we have that $\left\lfloor py \slash (4y - p) \right\rfloor + 1 = \left\lceil py \slash (4y-p) \right\rceil$.  Similarly for any prime $p \neq 2$  and $x \in \mathbb{N}$  that satisfy (\ref{eq: main})  we have that$\left\lfloor px \slash (4x - p) \right\rfloor + 1 = \left\lceil px \slash (4x-p) \right\rceil$.  This will help simplify how we express our work.  We now state our conjecture and provide a corollary to show the nature of our solution.

\begin{conjecture} \label{conj: gcd}
Consider a prime $p^{*} \geq 2521$.  Given any prime $p > p^{*}$  there exists $y \in \mathbb{N}$  so that $\left\lceil p \slash 2 \right\rceil \leq y \leq \left\lfloor p(p+3) \slash 6 \right\rfloor$, $\gcd \left( p,y \right) = 1$  and $$ \frac{y}{(4y-p)-m} \in \mathbb{N}$$
\ 

\noindent where $m \equiv py \mod (4y - p)$.
\end{conjecture}
\ 

\ 

\begin{corollary} \label{cor: gcd}
Consider a prime $p^{*} \geq 2521$.  Given any prime $p > p^{*}$  there exists $y \in \mathbb{N}$  so that $\left\lceil p \slash 2 \right\rceil \leq y \leq \left\lfloor p(p+3) \slash 6 \right\rfloor$, $\gcd \left( p,y \right) = 1$ and $$ \frac{4}{p} = \frac{1}{\left\lceil \frac{py}{4y-p} \right\rceil} + \frac{1}{y} + \frac{1}{p \cdot \lcm \left( \left\lceil \frac{py}{4y-p} \right\rceil ,y \right)}. $$
\end{corollary}
\ 

There are some scenarios for the prime $p$  that are guaranteed a solution of this type.  We outline the cases that have are guaranteed a solution in the following tables.  These results are incomplete and rather difficult to show in general.

\ 

\begin{center} 
\begin{tabular}{| c | c |} 
\hline 
p & y  \\ \hline 
$3 \mod 4$ & $(p(p+1) \slash 4)+1$ \\ \hline 
$5 \mod 8$ & $(3p+1) \slash 4$ \\ \hline 
$17 \mod 24$ & $(7p+1) \slash 4$ \\ \hline 
$97 \mod 120$ & $(7p+1) \slash 8$ \\ \hline 
$73 \mod 840$ & $(23p+1) \slash 8$ \\ \hline 
\end{tabular} 
\quad 
\begin{tabular}{| c | c |} 
\hline 
p & y  \\ \hline 
$241 \mod 840$ & $(23p+1) \slash 8$ \\ \hline 
$409 \mod 840$ & $(23p+1) \slash 8$ \\ \hline 
$433 \mod 840$ & $(15p+1) \slash 4$ \\ \hline 
$601 \mod 840$ & $(15p+1) \slash 4$ \\ \hline 
$769 \mod 840$ & $(15p+1) \slash 4$ \\ \hline 
\end{tabular}
\end{center}
\ 

\ 

We next make an analogue to conjecture \ref{conj: gcd} when the solutions are of type I(a).  This is much more enlightening for programming reasons.  We only need to check that the following conjecture holds for values of $x \in \mathbb{N}$  so that $\left\lceil p \slash 4 \right\rceil \leq x \leq \left\lfloor p \slash 2 \right\rfloor$.  The first conjecture will require us to search for a solution to (\ref{eq: main})  for values of $y$  on the boundary locations.  As $p$  gets large, the number of boundary locations grow at an asymptotic rate of $\mathcal{O}(p^2)$.  For this next conjecture, when considering type I(a) solutions, as $p$ gets large, the number of boundary locations grow at an asymptotic rate of $\mathcal{O}(p)$. This suggests that the result in \cite{sa} showing that every prime less than $10^{14}$  has a solution can be improved by searching for type I(a) solutions with $z = p \cdot \lcm(x,y)$.

Here we see that for all primes $p \not\in \{ 2, 3, 7, 47, 193, 2521 \}$  there exists a solution so that $y = \left\lceil px \slash (4x-p) \right\rceil$, $\gcd \left(p,y \right) = 1$  and $z = p \cdot \lcm \left( x,y \right)$.   We provide the foundation of this in the following conjecture. 

\begin{conjecture} \label{conj: x gcd}
Consider a prime $p^{*} \geq 2521$.  Given any prime $p > p^{*}$  there exists $x \in \mathbb{N}$  so that $\left\lceil p \slash 4 \right\rceil \leq x \leq \left\lfloor p \slash 2 \right\rfloor$, $\gcd \left( p, \left\lceil px \slash (4x-p) \right\rceil \right) = 1$  and $$ \frac{x}{(4x-p)-m} \in \mathbb{N}$$
\ 

\noindent where $m \equiv px \mod (4x - p)$.
\end{conjecture}

Much like conjecture \ref{conj: gcd}, this conjecture will lead to a solution of (\ref{eq: main}).  Now we will have the denominators of our unit fractions $x, \left\lceil px \slash (4x-p) \right\rceil$  and $p \cdot \lcm(x, \left\lceil px \slash (4x-p) \right\rceil )$.  We conclude our paper with more detail for some of our main points.  Section \ref{sec: development} is dedicated to some of the proofs to the propositions, corollaries and theorems made in this paper.

\section{Development} \label{sec: development}

\subsection{Proof of Proposition \ref{prop: three}}

\begin{proof}
\ 

\noindent Suppose that for a prime $p$  there exist values $x,y,z \in \mathbb{N}$  that make a solution to (\ref{eq: main}).  Further suppose that this solution is of type I(a). \\

\noindent This will imply that $$ y = \left\lfloor \frac{px}{4x-p} \right\rfloor +1.$$
\ 

\noindent We can clearly see that being a solution will imply that $$ \frac{4}{p} \geq \frac{1}{x} + \frac{1}{y}  $$
\ 

\noindent but to begin we will prove is that $$ \frac{4}{p} \leq \frac{1}{x-1} + \frac{1}{y}.$$

\noindent Proving this claim will lead us to show that it is a type I(b) solution. \\

\noindent First notice that for any prime $p$  and any $x \in \mathbb{N}$  such that $(p \slash 4) + 1 < x \leq (p \slash 2)$  we have that
\begin{align*}
\frac{p(x-1)}{4(x-1)-p} - \frac{px}{4x-p} &= \frac{p^{2}}{(4(x-1)-p)(4x-p))} \\
& \geq \frac{p^{2}}{(4x-p)^{2}} \\
& \geq 1.
\end{align*}
\ 

\noindent This tells us that 
\begin{align*} 
\frac{p(x-1)}{4(x-1)-p} &\geq \frac{px}{4x-p} + 1 \\
&\geq \left\lfloor \frac{px}{4x-p} \right\rfloor + 1.
\end{align*}
\ 

\noindent This will imply that 
\begin{align*} 
\frac{4}{p} &\leq \frac{1}{x-1} + \frac{1}{\left\lfloor \frac{px}{4x-p} \right\rfloor + 1} \\
&= \frac{1}{x-1} + \frac{1}{y}.
\end{align*}
\ 

\noindent To finish the proof we prove the following claim:  if $x,y,z \in \mathbb{N}$  is a solution to (\ref{eq: main}) for a prime $p$  and $$ \frac{4}{p} \leq \frac{1}{x-1} + \frac{1}{y} $$
\ 

\noindent then the solution is of type I(b). \\

\noindent  Because $$ \frac{4}{p} \geq \frac{1}{x} + \frac{1}{y}$$ 
\ 

\noindent and $$ \frac{4}{p} \leq \frac{1}{x-1} + \frac{1}{y}$$
\\ 

\noindent we see that $$ \frac{py}{4y-p} \leq x \leq \frac{py}{4y-p} + 1.$$
\ 

\noindent  Because $4xy - p(x+y) \neq 0$  for any $x,y \in \mathbb{N}$  that will make a solution to (\ref{eq: main}), we see that $py \slash (4y-p)$  is not an integer for the possible values of $y$  and $p$.  Because $x$  is a positive integer, we see then it must be true that $$ x = \left\lfloor \frac{py}{4y-p} \right\rfloor + 1.$$
\ 

\noindent This shows that the solution is of type I(b).
\end{proof}

\subsection{Proof of Theorem \ref{thm: two}}

\begin{proof}
\ 

\noindent This theorem is proved by the following selections of the value of $y$: \\

\centering 

\begin{tabular}{| c | c |}
\hline
p & y  \\ \hline
$2$ & $p(p+2) \slash 4$ \\ \hline
$3 \mod 4$ & $(p(p+1) \slash 4)+1$ \\ \hline
$5 \mod 8$ & $p(p+3) \slash 8$ \\ \hline
$17 \mod 24$ & $p(p+7) \slash 24$ \\ \hline
$73 \mod 120$ & $p(p+7) \slash 20$ \\ \hline
$97 \mod 120$ & $p(p+3) \slash 10$ \\ \hline
$4561 \mod 9240$ & $3p$ \\ \hline
$4729 \mod 9240$ & $3p$ \\ \hline
$5881 \mod 9240$ & $3p$ \\ \hline
$6049 \mod 9240$ & $3p$ \\ \hline
$6409 \mod 9240$ & $3p$ \\ \hline
$6841 \mod 9240$ & $3p$ \\ \hline
$7081 \mod 9240$ & $3p$ \\ \hline
$7729 \mod 9240$ & $3p$ \\ \hline
$8401 \mod 9240$ & $3p$ \\ \hline
$3049 \mod 9240$ & $p(p+31) \slash 44$ \\ \hline
$4369 \mod 9240$ & $p(p+31) \slash 44$ \\ \hline
$7009 \mod 9240$ & $p(p+31) \slash 44$ \\ \hline
$1201 \mod 9240$ & $5p(p+31) \slash 616$ \\ \hline
\end{tabular}
\quad
\begin{tabular}{| c | c |}
\hline
p & y  \\ \hline
$241 \mod 840$ & $p(p+11) \slash 42$ \\ \hline
$409 \mod 840$ & $p(p+11) \slash 42$ \\ \hline
$481 \mod 840$ & $p(p+11) \slash 84$ \\ \hline
$649 \mod 840$ & $p(p+11) \slash 84$ \\ \hline
$601 \mod 840$ & $p(p+15) \slash 56$ \\ \hline
$769 \mod 840$ & $p(p+15) \slash 56$ \\ \hline
$1009 \mod 9240$ & $3p$ \\ \hline
$1129 \mod 9240$ & $3p$ \\ \hline
$1801 \mod 9240$ & $3p$ \\ \hline
$2881 \mod 9240$ & $3p$ \\ \hline
$3649 \mod 9240$ & $3p$ \\ \hline
$4201 \mod 9240$ & $3p$ \\ \hline
$8521 \mod 9240$ & $3p$ \\ \hline
$8689 \mod 9240$ & $3p$ \\ \hline
$8929 \mod 9240$ & $3p$ \\ \hline
$3889 \mod 9240$ & $p(p+71) \slash 44$ \\ \hline
$5209 \mod 9240$ & $p(p+71) \slash 44$ \\ \hline
$7849 \mod 9240$ & $p(p+71) \slash 44$ \\ \hline
$6001 \mod 9240$ & $p(p+159) \slash 616$ \\ \hline
\end{tabular}
\ 

\ 

\flushleft

\noindent From this information you can derive the value of $z$  that solves equation \ref{eq: main}. \\

\ 

\noindent For example, if $p \equiv 5 \mod 8$, then there exists a value $k$  so that $p = 8k+5$.  We would see then that $y = (k+1)(8k+5)$. \\

\noindent Because 
\begin{align*}
\frac{py}{4y - p} &= \frac{(k+1)(8k+5)}{4k+3} \\
&= 2(k+1) - \frac{k+1}{4k+3} 
\end{align*}

\noindent and $0 < (k+1) \slash (4k+3) < 1$  for all $k \geq 0$, we see that $x = 2(k+1) = (p+3) \slash 4$. \\

\ 

\noindent Letting $x = (p+3) \slash 4$  and $y = p(p+3) \slash 8$  we see that necessarily $z = p(p+3) \slash 4$. \\

\ 

\noindent For every prime $p$  listed above, the given selection of $y$  will provide the values of $x$ and $z$ through the same process.
\end{proof}

\subsection{Proof of Corollary \ref{cor: gcd}}

\begin{proof}
\ 

\noindent If conjecture \ref{conj: gcd} holds then we necessarily have that $py \slash ((4y-p) - m) \in \mathbb{N}$  and one fact about every natural number $a \in \mathbb{N}$  is that $ \gcd (a, a+1) = 1$,  this will imply that $$ \gcd \left( \frac{py}{(4y-p)-m}, \frac{py}{(4y-p)-m} + 1 \right) = 1.$$
\ 

\noindent In particular, this would imply that $$ \gcd \left( py, py + (4y - p) - m \right) = (4y - p) - m.$$
\ 

\noindent Because $m \equiv py \mod (4y-p)$, we see that $$(4y-p) \left\lceil \frac{py}{4y-p} \right\rceil = py + (4y-p) - m.$$
\ 

\noindent This would imply that $$ \gcd \left( py, (4y-p) \left\lceil \frac{py}{4y-p} \right\rceil \right) = (4y-p) \left\lceil \frac{py}{4y-p} \right\rceil - py.$$
\ 

\noindent Because $\gcd (p,y) = 1$  we see that $\gcd ((4y-p),py) = 1$.  This will necessarily imply that  $$ \gcd \left( py, \left\lceil \frac{py}{4y-p} \right\rceil \right) = (4y-p) \left\lceil \frac{py}{4y-p} \right\rceil - py.$$
\ 

\noindent Because $ \left\lceil p \slash 4 \right\rceil \leq \left\lceil py \slash (4y-p) \right\rceil \leq \left\lfloor p \slash 2 \right\rfloor$  we see that $\gcd \left( \left\lceil py \slash (4y - p) \right\rceil, p \right) = 1$.  This will imply that $$ \gcd \left( y, \left\lceil \frac{py}{4y-p} \right\rceil \right) = (4y-p) \left\lceil \frac{py}{4y-p} \right\rceil - py.$$
\ 

\noindent We can express this as $$ 4y \left\lceil \frac{py}{4y-p} \right\rceil = py + p \left\lceil \frac{py}{4y-p} \right\rceil + \gcd \left( \left\lceil \frac{py}{4y-p} \right\rceil, y \right).$$
\ 

\noindent Dividing both sides of the equation by $py \left\lceil py \slash (4y - p) \right\rceil$, we have that $$ \frac{4}{p} = \frac{1}{\left\lceil \frac{py}{4y-p} \right\rceil} + \frac{1}{y} + \frac{1}{p \cdot \lcm \left( \left\lceil \frac{py}{4y-p} \right\rceil ,y \right)}. $$
\ 
\end{proof}

\bibliographystyle{amsplain}

\end{document}